\documentclass[10pt]{amsart}
\usepackage{latexsym}
\usepackage{amssymb}
\usepackage{amsmath}
\usepackage{mathrsfs}
\usepackage{bbm}
\usepackage{bbold}

\def\co{\operatorname{co}}

\def\spt{\operatorname{supp}}

\begin{document}

\title{
Pareto Optimality and Isoperimetry
}

\author{S.~S. Kutateladze}
\address[]{
Sobolev Institute of Mathematics\newline
\indent 4 Koptyug Avenue\newline
\indent Novosibirsk, 630090\newline
\indent Russia}
\email{
sskut@member.ams.org
}
\date{February 8, 2009}
%\sc{52A40; 90B50}
\begin{abstract}
Under study is the new class of geometrical extremal problems in which
it is required to achieve the best result in the presence of conflicting
goals; e.~g., given the surface area of a convex  body~$\mathfrak x$,
we try to maximize the volume of~$\mathfrak x$  and minimize the width
of ~$\mathfrak x$ simultaneously.
These problems are addressed along the lines of multiple criteria decision making.
\end{abstract}
\maketitle
%%%%%%%%%%%%%%%%%%%%%%%%%%%%%%%%%%%%%%%%%%%%%%%%%%%%%%%%%%%%%%%%%%%%%%%%%%

We address the multiple criteria
extremal problems of convex geometry which involve the goals and constraints
with the available description for the directional derivatives and the duals of the cones of feasible directions.
Transition to Pareto-optimality actually involves the scalar problems with
bulkier  objectives. The manner of combining the geometrical and functional-analytical
tools remains practically the same as in the case of a single goal typical of
an isoperimetric-type problem.
We  proceed by way of example and present here
a few model multiobjective problems that are connected
with the Blaschke and Minkowski structures.

Note that we use the notation and results of~\cite{MD}--\cite{Kut2007} as regards
convex geometry and the results of~\cite{Subdiff} as regards Pareto optimality.

{\bf 1.} {\scshape Vector Isoperimetric Problem}: Given are some convex bodies
$\mathfrak y_1,\dots,\mathfrak y_M$. Find a convex body $\mathfrak x$ encompassing a given volume
and minimizing each of the mixed volumes $V_1(\mathfrak x,\mathfrak y_1),\dots,V_1(\mathfrak x,\mathfrak y_M)$.
In symbols,
$$
\mathfrak x\in\mathscr A_N;\
\widehat p(\mathfrak x)\ge \widehat p(\bar{\mathfrak x});\
(\langle\mathfrak y_1,\mathfrak x\rangle,\dots,\langle\mathfrak y_M,\mathfrak x\rangle)\rightarrow\inf\!.
$$

Clearly, this is a~Slater regular convex program in the Blaschke structure.
Hence,  the following holds.

{\bf 2.} {\sl
Each Pareto-optimal solution $\bar{\mathfrak x}$ of the vector isoperimetric problem
has the form}
$$
\bar{\mathfrak x}=\alpha_1{\mathfrak y}_1+\dots+\alpha_m{\mathfrak y}_m,
$$
where $\alpha_1,\dots,\alpha_m$ are positive reals.

Let us illustrate 2 for the {\it Leidenfrost effect}, the
spheroidal state of a liquid drop on the horizontal heated surface.

{\bf 3.} {\scshape Leidenfrost Problem}. Given the volume of a three-dimensional
convex figure, minimize its  surface area and vertical  breadth.

By symmetry everything reduces to an analogous plane two-objective problem,
whose every Pareto-optimal solution is by~2 a~{\it stadium\/}, a weighted Minkowski sum of a disk and
a horizontal straight line segment.

{\bf 4.}
{\sl A plane spheroid, a Pareto-optimal solution of the Leidenfrost problem,
is the result of rotation of a stadium around the vertical axis through
the center of the stadium}.

{\bf 5.} {\scshape Internal Urysohn Problem with Flattening}.
Given are some~convex body
$\mathfrak x_0\in\mathscr V_N$ and some flattening direction~ $\bar z\in S_{N-1}$.
Among the convex bodies lying in~ $\mathfrak x_0$ and having
fixed integral breadth, find a convex body $\mathfrak x$
trying to maximize the volume of~$\mathfrak x$ and  minimize the
breadth of $\mathfrak x$ in the flattening direction:
$$
\mathfrak x\in\mathscr V_N;\
\mathfrak x\subset{\mathfrak x}_0;\
\langle \mathfrak x,{\mathfrak z}_N\rangle \ge \langle\bar{\mathfrak x},{\mathfrak z}_N\rangle;\
(-p(\mathfrak x), b_{\bar z}(\mathfrak x)) \to\inf\!.
$$

{\bf 6.}
{\sl For a feasible convex body $\bar{\mathfrak x}$ to be Pareto-optimal in
the internal Urysohn problem with the flattening
direction~$\bar z$ it is necessary and sufficient that there be
positive reals $\alpha, \beta$ and a~convex figure $\mathfrak x$ satisfying}
$$
\gathered
\mu(\bar{\mathfrak x})=\mu(\mathfrak x)+ \alpha\mu({\mathfrak z}_N)+\beta(\varepsilon_{\bar z}+\varepsilon_{-\bar z});\\
\bar{\mathfrak x}(z)={\mathfrak x}_0(z)\quad (z\in\spt(\mu(\mathfrak x)).
\endgathered
$$

By way of illustration we will derive the optimality criterion
in somewhat superfluous detail. In actuality, it would
suffice to appeal for instance to~\cite{Subdiff} or the other
numerous sources treating Pareto optimality in slightly less generality.

Note firstly that the internal Urysohn problem with flattening
may be rephrased in~ $C(S_{N-1})$ as the following
two-objective program
$$
\gathered
\mathfrak x\in \mathscr V_N;\\
\max\{{\mathfrak x}(z)- {\mathfrak x}_0(z)\mid z\in S_{N-1}\}\le 0;\\
\langle \mathfrak x,{\mathfrak z}_N\rangle \ge \langle\bar{\mathfrak x},{\mathfrak z}_N\rangle;\\
(-p(\mathfrak x), b_{\bar z}(\mathfrak x))\to\inf\!.
\endgathered
$$

The problem of Pareto optimization reduces to the
scalar program
$$
\gathered
\mathfrak x\in \mathscr V_N;\\
\max\{\max\{{\mathfrak x}(z)- {\mathfrak x}_0(z)\mid z\in S_{N-1}\},
\langle\bar{\mathfrak x},{\mathfrak z}_N\rangle-\langle \mathfrak x,{\mathfrak z}_N\rangle\}\le 0;\\
\max\{-p(\mathfrak x), b_{\bar z}(\mathfrak x)\}\to\inf\!.
\endgathered
$$

The last program is Slater-regular and so we may apply the
{\it Lagrange principle}. In other words, the value of the program under consideration
coincides
with the value of the unconstrained minimization problem
for an appropriate Lagrangian:
$$
\gathered
\mathfrak x\in \mathscr V_N;\\
\max\{-p(\mathfrak x), b_{\bar z}(\mathfrak x)\}\\ +\gamma\max\{
\max\{{\mathfrak x}(z)- {\mathfrak x}_0(z)\mid z\in S_{N-1}\},
\langle\bar{\mathfrak x},{\mathfrak z}_N\rangle-\langle \mathfrak x,{\mathfrak z}_N\rangle \}\to\inf\!.
\endgathered
$$
Here $\gamma$ is a~positive Lagrange multiplier.

We are left with differentiating the Lagrangian
along the feasible directions
and appealing to~the available results. Note in particular that
the relation
$
\bar{\mathfrak x}(z)={\mathfrak x}_0(z)\quad (z\in\spt(\mu(\mathfrak x))
$
is the  {\it complementary slackness condition\/}
 standard in mathematical programming.
The proof of the optimality criterion for the Urysohn problem with
flattening is complete.

Assume that a plane convex figure ${\mathfrak x}_0\in\mathscr V_2$ has the symmetry axis $A_{\bar z}$
with generator~$\bar z$.  Assume further that ${\mathfrak x}_{00}$ is the result of rotating
$\mathfrak x_0$
 around the symmetry axis $A_{\bar z}$ in~$\mathbb R^3$.  In this event we come to the following
 problem.

{\bf 7.} {\sc Internal Isoperimetric Problem in the class of
the surfaces of rotation with flattening in the direction of the axis of rotation}:
$$
\gathered
\mathfrak x\in\mathscr V_3;\\
\mathfrak x  \text{\ is\ a\ convex\ body\ of\ rotation\ around}\ A_{\bar z};\\
\mathfrak x\supset{\mathfrak x}_{00};\
\langle {\mathfrak z}_N, \mathfrak x\rangle \ge \langle{\mathfrak z}_N,\bar{\mathfrak x}\rangle;\\
(-p(\mathfrak x), b_{\bar z}(\mathfrak x)) \to\inf\!.
\endgathered
$$

By rotational symmetry, the three-dimensional problem reduces to an
analogous two-dimensional problem. The integral breadth and perimeter
are proportional on the plane, and we come to
the already settled problem~2. Thus, we have the following.

{\bf 8.}
{\sl Each Pareto-optimal solution of~7 is the result
of rotating around the symmetry axis a Pareto-optimal solution of the plane internal
Urysohn problem with flattening in the direction of the axis}.

Little is known about the analogous problems in arbitrary dimensions.
An especial place
is occupied by the result of Porogelov  who
demonstrated that the ``soap bubble'' in a tetrahedron
has the form of the result of the rolling of a ball over a~solution
of the internal Urysohn problem, i.~e. the weighted Blaschke sum of
a tetrahedron and a ball.

{\bf 9.} {\scshape External Urysohn Problem with Flattening}.  Given are some convex body
$\mathfrak x_0\in\mathscr V_N$ and some~flattening direction~$\bar z\in S_{N-1}$.
Among the
convex bodies encompassing ${\mathfrak x}_0$ and having fixed integral breadth,
find a convex body  $\mathfrak x$ maximizing value and minimizing  breadth in
the flattening direction:
$$
\mathfrak x\in\mathscr V_N;\
\mathfrak x\supset{\mathfrak x}_0;\
\langle \mathfrak x,{\mathfrak z}_N\rangle \ge \langle\bar{\mathfrak x},{\mathfrak z}_N\rangle;\
(-p(\mathfrak x), b_{\bar z}(\mathfrak x)) \to\inf\!.
$$

{\bf 10.}
{\sl For a feasible convex body $\bar{\mathfrak x}$ to be a Pareto-optimal solution of
the external Urysohn problem with flattening it is necessary and
sufficient that there be
positive reals $\alpha, \beta$, and a convex figure  $\mathfrak x$ satisfying}
$$
\gathered
\mu(\bar{\mathfrak x})+\mu(\mathfrak x)\gg{}_{{\mathbb R}^N} \alpha\mu({\mathfrak z}_N)+\beta(\varepsilon_{\bar z}+\varepsilon_{-\bar z});\\
V(\bar{\mathfrak x})+V_1(\mathfrak x,\bar{\mathfrak x})=\alpha V_1({\mathfrak z}_N,\bar{\mathfrak x})+ 2N\beta b_{\bar z}(\bar{\mathfrak x});\\
\bar{\mathfrak x}(z)={\mathfrak x}_0(z)\quad (z\in\spt(\mu(\mathfrak x)).
\endgathered
$$

Demonstration proceeds by analogy to the internal
Urysohn problem with flattening. The extra equality for mixed volumes
appears as the deciphering of the complementary slackness condition.

{\bf 11.} The above list  may be continued
with the multiobjective generalization of many scalar problems
such as problems with zone constraints and current hyperplanes,
problems over centrally symmetric convex figures,
Lindel\"of type problems, etc.  These problems
are usually convex with respect to Blaschke or Minkowski structures.
Of greater complexity
are the nonconvex parametric problems stemming from the
extremal properties of the Reuleaux triangle.
These problems require extra tools and undergone only a fragmentary
study.

In closing we dwell upon the problems of another type
where we seek for the  form of several convex figures simultaneously.

{\bf 12.} {\scshape Optimal Convex Hulls.}   Given are
convex bodies ${\mathfrak y}_1,\dots,{\mathfrak y}_m$ in~ $\mathbb R^N$.
Place a convex figure ${\mathfrak x}_k$ within~ ${\mathfrak y}_k$, for $k:=1,\dots,m$,
so as to simultaneously maximize the volume of
each of the figures $\mathfrak x_1,\dots,\mathfrak x_m$ and minimize
the integral breadth of the convex hull of the union of
these figures:
$$
\gathered
\mathfrak x_k\subset\mathfrak y_k \ (k:=1,\dots,m);\\
(-p({\mathfrak x}_1),\dots,-p({\mathfrak x}_m), \langle \co\{{\mathfrak x}_1,\dots,{\mathfrak x}_m\},{\mathfrak z}_N\rangle)\to\inf.
\endgathered
$$

{\bf13.}
{\sl For  some feasible convex bodies
${\bar{\mathfrak x}}_1,\dots,{\bar{\mathfrak x}}_m$  to have a
 Pareto-optimal convex hull it is necessary and sufficient
 that there be positive reals $\alpha_1,\dots,\alpha_m$ not
 vanishing simultaneously and two collections of
 positive Borel  measures $\mu_1,\dots,\mu_m$ and $\nu_1,\dots, \nu_m$ on~$S_{N-1}$ such that}
 $$
\gathered
\nu_1+\dots+\nu_m=\mu({\mathfrak z}_N);\\
 \bar{\mathfrak x}_k(z)={\mathfrak y}_k(z)\quad (z\in\spt(\mu_k));\quad\\
\alpha_k \mu(\bar{\mathfrak x}_k)=\mu_k+\nu_k\ (k:=1,\dots,m).
\endgathered
$$

The criterion appears along the lines
of 6.

\bibliographystyle{plain}

\enddocument